\documentclass[11pt]{amsart}

\usepackage{amssymb,amsthm,amsmath}
\usepackage[numbers,sort&compress]{natbib}
\usepackage{color}
\usepackage{graphicx}
\usepackage{tikz}
\usepackage[colorlinks,
linkcolor=red,
anchorcolor=blue,
citecolor=green
]{hyperref}

\newcommand{\SL}{\mathrm{SL}(2,\mathbb{R})}

\def\a{\alpha}

\def\e{\varepsilon}

\let\newpf\proof \let\proof\relax 
\newenvironment{pf}{\newpf[\proofname]}{\qed\endtrivlist}

\newcommand{\ba}{\overline{A}}

\def\be{\begin{equation}}
\def\ee{\end{equation}}

\def\ba{{\begin{align}}}
\def\ea{{\end{align}}}

\def\bm{\begin{matrix}}
\def\em{\end{matrix}}

\def\a{{\alpha}}

\def\SL{{\mathrm{SL}}}
\def\PSL{{\mathrm{PSL}}}

\def\0{{\mathbf 0}}

\newtheorem{Theorem}{Theorem}[section]
\newtheorem{Lemma}{Lemma}[section]
\newtheorem{Proposition}{Proposition}[section]

\newtheorem{Remark}{Remark}[section]

\newtheorem{Definition}{Definition}[section]

\numberwithin{equation}{section}

\theoremstyle{definition}

\newcommand{\dist}{\operatorname{dist}}

\newcommand{\C}{{\mathbb C}}

\newcommand{\Q}{{\mathbb Q}}
\newcommand{\R}{{\mathbb R}}
\newcommand{\T}{{\mathbb T}}

\newcommand{\Z}{{\mathbb Z}}

\def\B0{{\bold{0}}}

\catcode`\@=12

\def\Empty{}
\newcommand\oplabel[1]{
  \def\OpArg{#1} \ifx \OpArg\Empty {} \else
    \label{#1}
  \fi}

\newcommand{\comm}[1]{}
\newcommand{\comment}[1]{}

\begin{document}

\title[]{Stability of the non-critical spectral properties I:  arithmetic absolute continuity of the integrated density of states}

\author {Lingrui Ge}
\address{Department of Mathematics, University of California, Irvine CA, 92617
	} \email{lingruig@uci.edu}

\author {Svetlana Jitomirskaya}
\address{
Department of Mathematics, University of California, Irvine CA, 92617
} \email{szhitomi@uci.edu}

\author {Xin Zhao}
\address{Department of Mathematics, University of California, Irvine CA, 92617
	} \email{njuzhaox@126.com}

\begin{abstract}
We prove absolute continuity of the integrated density of states
for  frequency-independent analytic perturbations of the non-critical
almost Mathieu operator under arithmetic conditions on frequency.
\end{abstract}

\maketitle
\section{Introduction}
Analytic one-frequency Schr\"odinger operators on $\ell^2(\Z)$ are given by,
\begin{equation}\label{sch}
(H_{V,\alpha,x}u)_n=u_{n+1}+u_{n-1}+V(x+n\alpha)u_n,\ \ n\in\Z,
\end{equation}
where $\alpha\in\R\backslash\Q$ is the frequency, $x\in\T:=\R/\Z$ is
the phase, and $V\in C^\omega(\T,\R)$ \footnote{For a bounded
  1-periodic (possibly matrix valued) function $F$ analytic on
  $ \{ x |  | \Im x |< h \}$, and continuous on $ \{ x |  | \Im x
  |\leq  h \}$ let
$
\lvert F\rvert _h=  \sup_{ | \Im x |< h } \| F(x)\| $, and denote by $C^\omega_{h}(\T,*)$ the
Banach space of all these $*$-valued functions ($*$ will usually denote $\R$, $sl(2,\R)$,
$SL(2,\R)$).  $C^{\omega}(\T,\R)$ is the locally convex space given by
the union $\cup_{h>0}C_h^{\omega}(\T,\R)$.} is the potential.

The central and most extensively studied such operator is the almost Mathieu operator (AMO),
\begin{equation}\label{amo}
 (H_{\lambda,\alpha,x}u)_n=u_{n+1}+u_{n-1}+2\lambda\cos2\pi(x+n\alpha)u_n,\ \ n\in\Z.
\end{equation}
Sometimes called the drosophila of the subject. It is a model that is
responsible for both the origins of the field and much of its ongoing
significance in physics \cite{pei,hofst,sim,aos}. The almost Mathieu family is prototypical in
the sense of Avila's global theory, with its separate regimes
$\lambda<1,\;\lambda=1,$ and $\lambda>1$ lending names to global
classification of analytic $SL(2,\C)$ cocycles \cite{avila}.

At the
same time this family has a very special symmetry, self-duality with
respect to Aubry duality
(e.g. \cite{gjls}), that links $\{H_{\lambda,\alpha,x}\}_x$ and
$\{H_{\lambda^{-1},\alpha,x}\}_x$ and
stems from the gauge invariance of the underlying two-dimensional
model \cite{mz}. A number of remarkable results have been obtained, exploiting
this symmetry, thus with methods not extendable to the general
analytic class. This includes the ten martini problem
\cite{aj,ak,puig}, non-critical dry ten martini problem
\cite{aj1,ayz1}, the absolute continuity of the integrated density of
states in the non-critical case \cite{aj1,aviladamanik} and others.
See also \cite{avilaamo,jkr} for some other recent progresses. Aubry
duality enables one to combine the reducibility ($|\lambda|<1$) and
the localization ($|\lambda|>1$) techniques. However,  already for small perturbations of
the AMO,  most of the proofs involving self-duality of the AMO family
break down. It is therefore an interesting question whether self-duality is just
a convenient tool or of intrinsic importance to the AMO
results. Indeed, a few properties of the AMO are destroyed after
perturbations, for example the Lyapunov exponent is no longer a
constant on the spectrum, in general, making it particularly significant
 to identify those that
are stable.

 Another
(intersecting) group of important results exploits the fact that the potential
of the AMO is given by a first degree trigonometric  polynomial which
allows for some powerful considerations not available in the general
case. This includes metal-insulator transitions
\cite{J,ayz,JLiu,JLiu1,gyzh1}, and was also exploited in the ten
martini proofs \cite{puig, aj,aj1} and other arithmetic results, e.g.
\cite{gsv}. Here, by {\it arithmetic} we mean results under an explicit
arithmetic condition on the frequency.

In contrast, for general analytic one-frequency Schr\"odinger operators, the
current state-of-the-art results for all the above problems in the
positive Lyapunov exponent regime are measure-theoretic in $\alpha$
\cite{bg,gs2,gs3}. The biggest issue is that one needs to eliminate
frequencies $\alpha$ in a highly implicit way, technically due to the
need to get rid of the so-called
``double resonances''.

At the same time, the almost Mathieu operator coming from physics, it is
natural to expect that its physically relevant properties hold at
least for its small perturbations. In this respect it is particularly
important that the allowed perturbations are uniform in $\alpha$
within the Diophantine class \footnote{A result where the strength of
  the allowed perturbation
depends, say, on the Diophantine constants of $\alpha$
is clearly not robust with respect to small changes of the Diophantine
frequency and requires various positive measure exclusions if
a.e. frequency is fixed.}. Results with uniform dependence on
Diophantine $\alpha$ are often called {\it non-perturbative}
(e.g. \cite{bj02}), even when the parameters involved are otherwise small, while the ones
without such dependence are called {\it perturbative}.

The paper is the first of a multi-part project to extend various
spectral properties of the almost Mathieu operator that have so far been proved in an
AMO-specific way, to the analytic neighborhood of the almost
Mathieu operators in a {\it
  non-perturbative} way, so that, in particular, to confirm their relevance to physics.
Namely, we consider
\begin{equation}\label{pamo}
(H^\e_{\lambda,\alpha,x}u)_n=u_{n+1}+u_{n-1}+(2\lambda\cos2\pi(x+n\alpha)+\e v(x+n\alpha))u_n,\ \ n\in\Z,
\end{equation}
where $v$ is a 1-periodic real analytic function. An additional aim is to
develop different techniques to investigate the  spectral properties
of operators \eqref{pamo} in the zero/positive Lyapunov exponent regimes in a way that
does not use self-duality or low degree of the potential, with an
expectation that some of the techniques will also turn out to be
useful for global results.

Technically, from the point of view of Avila's global theory \cite{avila}, the AMO
family has one more important feature: the acceleration (see \eqref{ac})
is bounded by $1$ on the spectrum.
An important goal of the present project is to show that it is exactly
this feature that governs many of the spectral properties for
operators \eqref{pamo}, that prevents, in particular, the occurrence of
the double resonances, thus confirming the importance of the notion of
acceleration in the spectral theory
of analytic quasiperiodic operators.

The {\it acceleration} is defined as
\begin{equation}\label{ac}
\omega(E)=\lim\limits_{\epsilon\rightarrow 0^+}\frac{L_\epsilon(E)-L_0(E)}{2\pi\epsilon},
\end{equation}
where $L_\epsilon(E)$ is the complexified Lyapunov exponent:
\begin{align}\label{multiergodicsch}
L_\epsilon(E)=\lim\limits_{n\rightarrow\infty}\frac{1}{n} \int_\T \ln \|A(x+i\epsilon+(n-1)\alpha)\cdots A(x+i\epsilon)\|dx,
\end{align}
with
$$
A(x)=\begin{pmatrix}E-V(x)&-1\\ 1&0\end{pmatrix}.
$$
One of the key conclusions in \cite{avila} is that the acceleration is always an integer. Moreover, for the almost Mathieu operator, for all $E$ in the spectrum, we have
\begin{enumerate}
\item $|\lambda|<1$: $L(E)=0$ and $\omega(E)=0$.
\item $|\lambda|=1$: $L(E)=0$ and $\omega(E)=1$.
\item $|\lambda|>1$: $L(E)=\ln|\lambda|>0$ and $\omega(E)=1$.
\end{enumerate}
For general one-frequency Schr\"odinger operators, one can similarly divide  the spectrum into three regimes:
\begin{enumerate}
\item The subcritical regime: $L(E)=0$ and $\omega(E)=0$.
\item The critical regime: $L(E)=0$ and $\omega(E)>0$.
\item The supercritical regime: $L(E)>0$ and $\omega(E)>0$.
\end{enumerate}

Now, for  the analytic perturbations of the non-critical almost
Mathieu operator \eqref{pamo}, as was also proved in \cite{avila}, for all $E$ in the spectrum, one has
\begin{enumerate}
\item $|\lambda|<1$ and  $\e$ small enough: $L(E)=0$ and $\omega(E)=0$.
\item $|\lambda|>1$ and $\e$ small enough: $L(E)>0$ and $\omega(E)=1$.
\end{enumerate}

 We would like to mention that the spectral properties of non-critical
operators \eqref{pamo}  were well studied in the perturbative regime, i.e., assuming
 $|\lambda|$ sufficiently large depending on $\alpha$. For a {\it
   fixed} Diophantine frequency, we refer readers to
 \cite{Sinai,fsw,fv} for the proofs of almost sure Anderson
 localization, to  \cite{gyzh2} for the arithmetic version of Anderson
 localization in the case the potential is even, to \cite{wz1,fv} for
 the proofs of Cantor spectrum, and to \cite{B,wz,LWY,xgw} for the
 positivity and (H\"older) continuity of the Lyapunov exponent, all
 even holding for much rougher $C^2$ potentials as long as they
 stay $\cos$-type.

In the present paper, we want to study  the regularity of the
integrated density of states (IDS) in the global sense ($|\lambda|$
does not need to be large) and non-perturbatively.

The IDS is defined in a uniform way for the above one-frequency Schr\"odinger operators $(H_{V,\alpha,x})_{x\in\T}$ by
$$
N(E)=\int_{\T}\mu_{x}(-\infty,E]dx,
$$
where $\mu_x$ is the spectral measure associated with $H_{V,\alpha,x}$ and $\delta_0$. Roughly speaking, the density of states measure $N([E_1,E_2])$ gives the ``number of states per unit volume" with energy between $E_1$
and $E_2$.

Regularity of the IDS is a popular subject in the spectral theory
of quasiperiodic operators, especially the absolute continuity
\cite{aviladamanik,avila1,aj1} and the H\"older regularity
\cite{amor,aj1,gs1,gs2}. It is also closely related to many other
topics. For example, absolute continuity of the  IDS is closely
related to purely absolutely continuous spectrum in the regime of zero
Lyapunov exponent \cite{kotani,damanik}. H\"older continuity of the IDS is
closely related to homogeneity of the spectrum
\cite{dgl,dgsv,lyzz}. Before formulating our results, we first give
precise arithmetic assumptions on $\alpha$. A frequency  $\alpha\in
\R$ will be called ($\kappa,\tau$)-{\it strongly Diophantine}   (denoted by $\alpha\in \rm{SDC}(\kappa,\tau))$ where $\kappa>0$, $\tau>1$ if
\begin{equation}\label{dc_def}
\dist(k\alpha,\Z)\ge \frac{\kappa}{|k|(\ln|k|)^{\tau}},\quad \forall k\in \Z\backslash\{0\}.
\end{equation}
We will use the notation
$$
\rm{SDC}:=\bigcup_{\kappa>0;\,\tau>1} \rm{SDC}(\kappa,\tau).
$$
Clearly, $\rm{SDC}$ is a set of  full Lebesgue measure.

We have
\begin{Theorem}\label{main}
Let $\alpha\in {\rm SDC}$, $|\lambda|\neq 1$ and $v$ be real analytic.  There is $\e_0(\lambda,v)>0$, such that if $|\e|<\e_0$, then the integrated density of states of operator \eqref{pamo} is absolutely continuous.
\end{Theorem}

\begin{Remark}
The strong Diophantine condition on $\alpha$ can be relaxed to the
usual Diophantine condition, or even to the Bruno condition where we
believe the method in this paper still works. We only use it to invoke
the existing homogeneity results and to keep
the paper short.
\end{Remark}
\begin{Remark}
This paper mainly deals with the perturbations of the supercritical
AMO. For the perturbations of the subcritical AMO, absolute continuity
of the IDS is a corollary of the almost reducibility conjecture (ARC),
announced by Avila \cite{avila}. However, to keep the paper
self-contained, we give a short proof of it, independent of the ARC.
\end{Remark}
\begin{Remark}
If $\lambda=0$, Theorem \ref{main} follows from \cite{aj1}.
\end{Remark}
\begin{Remark}
Note that $\e_0$ actually depends  only on $\lambda$ and the analytic norm of $v$.
\end{Remark}

In fact, we only use a single special feature of the AMO, and effectively prove the following theorem:

\begin{Theorem}\label{main0}
Suppose every $E$ in the spectrum of operator $H_{V,\alpha,x}$ with
real-analytic $V$ given by
\eqref{sch} is non-critical, and satisfies $\omega(E)\leq 1.$ Let $W$ be real analytic.  There is
$\e_0(V,W)>0$, such that if $\alpha\in
{\rm SDC}$ and $|\e|<\e_0$, then the integrated
density of states of operator $H_{V+\e  W,\alpha,x}$ is absolutely continuous.
\end{Theorem}

\begin{Remark}
 Proof of Theorem \ref{main} works verbatim to prove the
  supercritical part of Theorem \ref{main0}. As for the subcritical
  part of Theorem \ref{main0}, one needs to invoke the ARC, because
  our short self-contained argument in Section \ref{sub} only works
  for perturbations of the
  AMO.
\end{Remark}

Finally, we briefly introduce the main ideas of our proof.  Previously
methods to prove
absolute continuity
of the IDS in the supercritical regime were developed in \cite{gs2} using  the large deviation theorem
and avalanche principle for determinants of the truncated operator. However, this only works in a
measure-theoretic sense (one cannot fix $\alpha$). In the subcritical
regime, the absolute continuity of the IDS is, as mentioned, a corollary of the
almost reducibility. Here, the method we use is rather
different. For operators with homogenous spectrum, Sodin and Yuditskii
\cite{sy1,sy2} gave a characterization of the absolute continuity of
the corresponding spectral measure: it is purely absolutely continuous
if the real part of the normal boundary of the Borel transform of the
spectral measure is integrable and the topological support of it is
homogeneous. This characterization can be verified if the operator is
reflectionless since the real part of the normal boundary is zero. In
this paper, we use Avila's global theory \cite{avila} to develop a
new way to study the real part of the normal boundary of the Borel
transform of the IDS, which is enough to show that it is integrable.

\section{Preliminaries}

\subsection{Quasiperiodic cocycles and the Lyapunov exponent}

Given $A \in C^\omega(\T,{\rm SL}(2,\R))$ and $\alpha\in\R\backslash\Q$, we define the  \textit{quasiperiodic $SL(2,\R)$-cocycle} $(\alpha,A)$:
$$
(\alpha,A)\colon \left\{
\begin{array}{rcl}
\T \times \R^2 &\to& \T \times \R^2\\[1mm]
(x,v) &\mapsto& (x+\alpha,A(x)\cdot v)
\end{array}
\right.  .
$$
The iterates of $(\alpha,A)$ are of the form $(\alpha,A)^n=(n\alpha,  \mathcal{A}_n)$, where
$$
\mathcal{A}_n(x):=
\left\{\begin{array}{l l}
A(x+(n-1)\alpha) \cdots A(x+\alpha) A(x),  & n\geq 0\\[1mm]
A^{-1}(x+n\alpha) A^{-1}(x+(n+1)\alpha) \cdots A^{-1}(x-\alpha), & n <0
\end{array}\right.    .
$$
The {\it Lyapunov exponent} is defined by
$\displaystyle
L(\alpha,A):=\lim\limits_{n\to \infty} \frac{1}{n} \int_{\T} \ln \|\mathcal{A}_n(x)\| dx.
$

A basic fact about  quasiperiodic $SL(2,\R)$-cocycle is the continuity of the Lyapunov exponent:
\begin{Theorem}[\cite{bj}]\label{lyacon}
The functions $\R \times C^{\omega}(\T, {\rm SL}(2,\R))\ni
(\alpha,A)\mapsto L(\alpha,A)\in [0,\infty)$
are continuous at any $(\alpha',A')$ with $\alpha'\in \R\backslash\Q$.
\end{Theorem}
\subsection{The rotation number}
Assume that $A \in C^\omega(\T, {\rm SL}(2, \R))$ is homotopic to the identity. $(\alpha, A)$ induces the projective skew-product $F_A\colon \T \times \mathbb{S}^1 \to \T \times \mathbb{S}^1$
$$
F_A(x,w):=\left(x+\a,\, \frac{A(x) \cdot w}{|A(x) \cdot w|}\right),
$$
which is also homotopic to the identity.
Thus we can lift $F_A$ to a map $\widetilde{F}_A\colon \T \times \R \to \T \times \R$ of the form $\widetilde{F}_A(x,y)=(x+\alpha,y+\psi_x(y))$, where for every $x \in \T$, $\psi_x$ is $\Z$-periodic.
The map $\psi\colon\T \times \T  \to \R$ is called a {\it lift} of
$A$. Let $\mu$ be any probability measure on $\T \times \R$ which is
invariant by $\widetilde{F}_A$, and whose projection on the first
coordinate is given by the Lebesgue measure.
The number
$$
\rho(\alpha,A):=\int_{\T \times \R} \psi_x(y)\ d\mu(x,y) \ {\rm mod} \ \Z
$$
 depends  neither on the lift $\psi$ nor on the measure $\mu$, and is called the \textit{fibered rotation number} of $(\alpha,A)$ (see \cite{H,JM} for more details).

A typical  example is represented by the \textit{Schr\"{o}dinger cocycles} $(\alpha,S_E^{V})$, where
$$
S_E^{V}(x):=
\begin{pmatrix}
E-V(x) & -1\\
1 & 0
\end{pmatrix},   \quad E\in\R.
$$
Schr\"odinger cocycles are a dynamical equivalent of the eigenvalue equations $H_{V, \alpha, x}u=E u$. Indeed, any formal solution $u=(u_n)_{n \in \Z}$ of $H_{V, \alpha, x}u=E u$ satisfies
$$
\begin{pmatrix}
u_{n+1}\\
u_n
\end{pmatrix}
= S_E^V(x+n\alpha) \begin{pmatrix}
u_{n}\\
u_{n-1}
\end{pmatrix},\quad \forall \  n \in \Z.
$$
The spectral properties of $H_{V,\alpha,x}$ and the dynamics of
$(\alpha,S_E^V)$ are closely related by the Johnson's theorem:
 $E\in \Sigma_{V,\alpha}$ \footnote{
By minimality, the spectrum of $H_{V,\alpha,x}$ denoted by
$\Sigma_{V,\alpha}$, is a compact subset of $\R$, independent of $x$
if $(1,\alpha)$ is rationally independent.}  if and only if
$(\alpha,S_E^{V})$ is \textit{not} uniformly hyperbolic. Throughout
the rest of the paper, we set $L(E):=L(\alpha,S_E^{V})$  \footnote{We sometimes identify $L_0(E)$ (see \eqref{multiergodicsch}) and $L(E)$.}  and $\rho(E):=\rho(\alpha,S_E^{V})$ for brevity.  It is also well known that $\rho(E)\in[0,\frac{1}{2}]$ relates to the integrated density of states as follows:
\begin{equation}\label{relation}
N(E)=1-2\rho(E).
\end{equation}

\section{Proof for the supercritical case.}
The proof of Theorem \ref{main} contains two parts, based on two different methods. In this section, we deal with the perturbations of the supercritical almost Mathieu operators.
\begin{Theorem}\label{sup}
Let $\alpha\in {\rm SDC}$, $|\lambda|>1$ and $v$ be real analytic.  Then for $\e$ small enough, depending on $\lambda,v$, the IDS of operator \eqref{pamo} is absolutely continuous.
\end{Theorem}
Define
$$
\delta_j(n)=\begin{cases}
1& n=j\\
0& n\neq j
\end{cases}.
$$
For any $E\in\R$ and $\epsilon>0$, one can define the averaged Green's function of operator \eqref{sch} as
$$
G(0,0,E+i\epsilon)=\int_\T \langle\delta_0, (H_{V,\alpha,x}-(E+i\epsilon))^{-1}\delta_0\rangle dx=\int \frac{1}{E'-(E+i\epsilon)}dN(E').
$$
Let us recall that a {\it Herglotz function} is a holomorphic mapping of $\C^+=\{z\in\C:\Im z>0\}$ to itself. One can easily check that $G(0,0,z)$ is a Herglotz function. Thus for almost every $E\in\R$, the normal boundary of $G(0,0,z)$ exists and one can define
$$
G(0,0,E+i0)=\lim_{\epsilon\rightarrow 0^+}G(0,0,E+i\epsilon).
$$
Before we give the proof of Theorem \ref{sup}, let's recall two interesting theorems.  Given a compact set $\mathcal{S}$, we say $\mathcal{S}$ is homogenous if there is $\sigma_0>0$ such that for any $\sigma<\sigma_0$ and $E\in \mathcal{S}$, we have
$$
\left|(E-\sigma,E+\sigma)\cap\mathcal{S}\right|\geq \frac{1}{2}\sigma.
$$
\begin{Theorem}[Theorem H of \cite{dgsv}]\label{homo}
Assume $\alpha\in {\rm SDC}$, $V$ is real analytic and $L(E)>0$ for all $E\in\R$, then $\Sigma_{V,\alpha}$ is homogeneous.
\end{Theorem}
\begin{Remark}
As explained in \cite{dgsv}, the strong Diophantine condition on
$\alpha$ can be relaxed to the usual Diophantine condition.
\end{Remark}
\begin{Theorem}[Lemma 2.4 of \cite{sy1}]\label{abs}
Let $\mathcal{E}\subset \R$ be a compact homogenous set and $f$ a Herglotz function with representation
$$
f(z)=\int_{\mathcal{E}}\frac{d\mu(E)}{E-z}, \ \ z\in\C_+,
$$
where $d\mu$ is a finite measure with $supp(d\mu)\subset
\mathcal{E}$. Let $f(E+i0)=\lim_{\epsilon\rightarrow
  0^+}f(E+i\epsilon)$ be the a.e. normal boundary of $f$ and assume that
$$
\Re f(E+i0)\in L^1(\mathcal{E},dE).
$$
Then $d\mu$ is absolutely continuous.
\end{Theorem}
Theorem \ref{homo} and Theorem \ref{abs} indicate that to prove absolute continuity of the IDS, we need to study the regularity of the real part of the normal boundary of the averaged Green's function.
It's easy to see that $\Re G(0,0,E+i0)$ is well-defined and real
analytic outside the spectrum. Beyond that, the
key statement of the celebrated Kotani theory \cite{kotani}, says
that $ G(0,0,z)$ is {\it reflectionless}  in the zero Lyapunov exponent
regime. This means that for almost every $E$ in the zero Lyapunov
exponent regime, $\Re G(0,0,E+i0)=0$.

However, in the positive
Lyapunov exponent regime, the regularity of $\Re G(0,0,E+i0)$ remains
widely open. Indeed, one can expect that in this case $\Re
G(0,0,E+i0)$ as a function of $E$ can be as bad  as possible, since
the spectrum is purely singular. It turns out that this standard
intuition is completely wrong. Indeed, we have the following surprising
theorem.
\begin{Theorem}\label{main1}
Let $\alpha\in\R\backslash\Q$, $|\lambda|\neq 0,1$ and $v$ be real analytic.  Then for $\e$ small enough, $\Re G_\e(0,0,E+i0)$ \footnote{$$
G_\e(0,0,E+i0)=\lim\limits_{\epsilon\rightarrow 0^+}G_\e(0,0,E+i\epsilon),
$$
and
$$
G_\e(0,0,E+i\epsilon)=\int_\T \langle\delta_0, (H^\e_{\lambda,\alpha,x}-(E+i\epsilon))^{-1}\delta_0\rangle dx.
$$} almost surely coincides with an analytic function on the spectrum.
\end{Theorem}
Here we say $f$ is a real analytic function on a set $S$ if it is
the restriction of some real analytic function defined on an open
neighborhood of $S.$
\begin{Remark}
Reflectionless property is just the same result, but with the analytic function
being $0$. In this sense,  almost sure  analyticity can be viewed as
the generalized notion of {\it reflectionless}.
\end{Remark}
Note that Theorem \ref{sup} follows from Theorem \ref{main1}.\\
{\bf Proof of Theorem \ref{sup}:} For the supercritical AMO, it is
well known that the Lyapunov exponent is positive for all $E\in\R$ \cite{bj,avila}. Thus by Theorem \ref{lyacon}, there is $\e_0(\lambda,v)>0$, such that if $|\e|<\e_0$, then
$$
L(\alpha, S_E^{2\lambda\cos+\e v})>0
$$
for any $E\in\R$. By Theorem \ref{homo}, if $\alpha\in{\rm SDC}$, then $\Sigma^\e_{\lambda,\alpha}$ is homogenous where $\Sigma^\e_{\lambda,\alpha}$ is the spectrum of $H_{\lambda,\alpha,x}^\e$.  On the other hand, by Theorem \ref{main1}, there is $\e_1(\lambda,v)>0$, such that if $|\e|<\e_1$, then $\Re G_\e(0,0,E+i0)$ almost surely equals to an analytic function on the spectrum. Thus for $|\e|<\min\{\e_0,\e_1\}$, we have $\Sigma_{\lambda,\alpha}^\e$ is homogenous and $\Re G_\e(0,0,E+i0)\in L^1(\Sigma_{\lambda,\alpha}^\e)$. Note that
$$
G_\e(0,0,z)=\int\frac{1}{E-z}dN_\e(E), \ \ z\in\C_+.
$$
Thus by Theorem \ref{abs}, we have $N_\e(E)$, the IDS of operator \eqref{pamo}, is absolutely continuous.\qed

In the remaining part of this section we prove Theorem
\ref{main1}. The foundation  is the following remarkable result of Avila in \cite{avila}, on the analyticity of the Lyapunov exponent.
\begin{Theorem}[\cite{avila}]\label{sle}
Let $\lambda>1$ and $v$ be any real analytic function. Then for $\e$ small enough and for every $\alpha\in\R\backslash\Q$, $L(\alpha,S_E^{2\lambda \cos+\e v})$ restricted to the spectrum is a positive real analytic function.
\end{Theorem}
On the other hand, there is the following relation between the
Lyapunov exponent and the Green's function, see \cite{JM,kotani,ks},
$$
\frac{\partial L(E+i\epsilon)}{\partial E}=\Re G(0,0,E+i\epsilon).
$$
Roughly speaking, the real part of the normal boundary of the averaged
green's function is exactly the normal boundary of the derivative of
the Lyapunov exponent. We will now show that $\frac{\partial L(E+i0)}{\partial E}$ is almost surely analytic. To prove this, we will need some  ideas from hard analysis.

\subsection{Non-tangential maximal functions}
\begin{Definition}
The non-tangential maximal function takes a function $F$ defined on the upper-half plane
$$
\C_+:=\{x+iy:x\in\R, y>0\}
$$
and produces a function $F^*$ defined on $\R$ via the expression
$$
F^*(x_0)=\sup\limits_{|x-x_0|<y} |F(x+iy)|.
$$
\end{Definition}
Note that for any fixed $x_0$, the set $\{(x,y):|x-x_0|<y\}$ is a cone in $\R^2_+$ with vertex at $(x_0,0)$ and axis perpendicular to the boundary of the $x$-axis. Thus, the non-tangential maximal operator simply takes the supremum of the function $F$ over a cone with vertex at the $x$-axis.
\begin{Definition}
The Hardy space $H^p$ where $0<p<\infty$, on the upper half-plane $\C_+$ is defined to be the space of holomorphic functions $F$ on $\C_+$ with bounded norm, the norm being given by
$$
|F|_{H^p}=\sup\limits_{y>0}\left(\int_{-\infty}^{+\infty}|F(x+iy)|^pdx\right)^{\frac{1}{p}}.
$$
\end{Definition}
\begin{Proposition}[Page 1 of  \cite{hl}, Theorem 1 of \cite{bgshar}]\label{hl}
Let  $F$ be an analytic function on the upper-half plane, and of the class $H^p$ where $0<p<\infty$. Then
$$
\int_{-\infty}^\infty(F^*(x))^pdx\leq C_p \sup\limits_{y>0}\int_{-\infty}^{+\infty}|F(x+iy)|^pdx.
$$
\end{Proposition}
Note that $G(0,0,z)$ is an analytic function on the upper-half plane. Thus one can define the corresponding non-tangential maximal function
$$
G^*(E)=\sup\limits_{|E'-E|<\epsilon} |G(0,0,E'+i\epsilon)|.
$$
\begin{Proposition}\label{em}
For each $\sigma>0$,
$$
\left|\{E: G^*(E)>\sigma\}\right|\leq \frac{D}{\sigma^{\frac{3}{4}}},
$$
for some $D>0$ (does not depend on $\sigma$).
\end{Proposition}
\begin{pf}
One can check that $G(0,0,\cdot)\in H^{p}$ for any $\frac{1}{2}<p<1$. Thus by Proposition \ref{hl}, there is $D>0$ such that
$$
\int_{-\infty}^\infty(G^*(E))^\frac{3}{4}dE\leq D.
$$
Thus
\begin{align*}
\sigma^{\frac{3}{4}} \left|\{E: G^*(E)>\sigma\}\right|&\leq \int_{\{E: G^*(E)>\sigma\}} (G^*(E))^{\frac{3}{4}}dE\\
&\leq \int_{-\infty}^\infty(G^*(E))^\frac{3}{4}dE\leq D.
\end{align*}
\end{pf}

\subsection{Proof of Theorem \ref{main1}} Now we are ready to prove Theorem \ref{main1}. For simplicity, we will omit $\e$ in the notations. Note that by the spectral theorem, for any $z\in \C\backslash \R$, we have
$$
G(0,0,z)=\int \frac{1}{E-z} dN(E).
$$
On the other hand, denoting
$$
w(z)=\int \ln (E-z) d N(E).
$$
The Thouless formula \cite{craigsimon1,craigsimon2,haro} says
$$
L(z)=\int \ln|E-z|dN(E)=\Re w(z).
$$
The followings are some basic facts on $G^*(E)$,
\begin{enumerate}
\item $G^*$ is lower semicontinuous and in particular,
$$
U_\sigma=\{E: G^*(E)>\sigma\}
$$
is open.
\item By Proposition \ref{em}, $\R\doteq\cup_n U^c_n$ where we say $A\doteq B$ if $|A\backslash B|=|B\backslash A|=0$.
\end{enumerate}
Note that $C_n:=U_n^c\cap [-n,n]$ is compact. We define
$$
\epsilon(E)={\rm dist}(E, C_n).
$$
One can easily verify that
$$
|\epsilon(E)-\epsilon(E')|\leq |E-E'|
$$
for any $E,E'\in\R$. By  Rademacher's theorem, $\epsilon(E)$ is an absolutely continuous function and is differentiable almost surely.

On the other hand, $C_n^c$ is open, thus $C^c_n=\cup_i (a_i,b_i)$, and
we have
\begin{align}
\epsilon'(E)=\left\{
\begin{aligned} &0 &E\in C_n\backslash \left(\{a_i\}\cup\{b_i\}\right),\\
&1
& E\in \cup_i \left(a_i,\frac{a_i+b_i}{2}\right),\\
&-1
& E\in \cup_i \left(\frac{a_i+b_i}{2},b_i\right).
\end{aligned}\right.
\end{align}

For any fixed $\delta>0$, we consider the function
$f_\delta(E)=L\left(E+i(\epsilon(E)+\delta)\right)$, which is obviously Lipschitz. Thus
$$
f_\delta(E_1)-f_\delta(E_2)=\int_{E_1}^{E_2}\frac{\partial L}{\partial E} \left(E+i(\epsilon(E)+\delta)\right)+\frac{\partial L}{\partial \epsilon} \left(E+i(\epsilon(E)+\delta)\right) \epsilon'(E)dE.
$$
The following Lemma plays a crucial role.
\begin{Lemma}\label{l51}
For any $E_1,E_2\in C_n$, we have
$$
L(E_1)-L(E_2)=\int_{E_1}^{E_2} g(E)dE,
$$
where
\begin{equation}\label{forg}
g(E)=\begin{cases}
\lim\limits_{\delta\rightarrow 0^+} \frac{\partial L}{\partial E}(E+i\delta) & E\in C_n,\\
\frac{\partial L}{\partial E} \left(E+i\epsilon(E)\right)+\frac{\partial L}{\partial \epsilon} \left(E+i\epsilon(E)\right) \epsilon'(E) &  E\notin C_n.
\end{cases}
\end{equation}
\end{Lemma}
\begin{pf}
Since $E_1,E_2\in C_n$, for all $E\in (E_1,E_2)$ and all $\delta>0$ small, by the definition of $\epsilon(E)$, we have
$$
(E,\epsilon(E)+\delta)\in \bigcup\limits_{E\in C_n} \left\{(E',\epsilon'): |E'-E|\leq \epsilon'\right\}.
$$
Thus by the definition of $G^*(E),C_n$ and the fact that
\begin{equation*}
\frac{dL}{dz}(E+i\epsilon)=G(0,0,E+i\epsilon),
\end{equation*}
we have
$$
\left|\frac{\partial L}{\partial E} \left(E+i(\epsilon(E)+\delta)\right)+\frac{\partial L}{\partial \epsilon} \left(E+i(\epsilon(E)+\delta)\right) \epsilon'(E)\right|\leq 2n,
$$
uniformly for all $\delta>0$. Thus the result follows from dominated convergence by letting $\delta\rightarrow 0$.
\end{pf}
Now, we apply Lebesgue's theorem on differentiation of integrals to $g\in L^\infty$,
\begin{equation}\label{s1}
g(E)=\lim\limits_{E'\rightarrow E}\frac{1}{E'-E}\int_E^{E'} g(E)dE,
\end{equation}
for a.e. $E\in\R$.

We may assume $|C_n\cap \Sigma^\e_{\lambda,\alpha}|>0$ (otherwise there is nothing to say). Then by Lebesgue's density theorem,  for a.e. $E\in C_n\cap\Sigma^\e_{\lambda,\alpha}$, there is a sequence $E_j\in C_n\cap\Sigma^\e_{\lambda,\alpha}$ such that
\begin{equation}\label{s2}
\lim_{j\rightarrow \infty} E_j=E.
\end{equation}

Now combining \eqref{s1}, \eqref{s2} with Lemma \ref{l51}, for a.e. $E\in C_n\cap\Sigma_{\lambda,\alpha}^\e$, we can find a sequence $E_j\in C_n\cap\Sigma_{\lambda,\alpha}^\e$ such that
\begin{equation}\label{weneed}
g(E)=\lim\limits_{j\rightarrow \infty}\frac{L(E_j)-L(E)}{E_j-E}.
\end{equation}
Note that by Theorem \ref{sle}, there is an analytic function $f:U\rightarrow \R$ such that $L(E)=f(E)$ on $\Sigma_{\lambda,\alpha}^\e\subset U$ where $U$ is open. Thus by \eqref{forg} and \eqref{weneed}, we have for a.e. $E\in C_n\cap\Sigma_{\lambda,\alpha}^\e$,
$$
\Re G(0,0,E+i0)=\lim\limits_{\epsilon\rightarrow 0^+} \frac{\partial L}{\partial E}(E+i\epsilon)=g(E)=f'(E).
$$
Note that $f'(E)$ is also analytic. This completes the proof.

\section{Proof of the case $|\lambda|<1$}\label{sub}
Absolute continuity of the IDS in the
subcritical regime is actually  a corollary of the almost
reducibility \cite{avila1,avila2}. Thus the main aim of the section is
to prove that the cocycle corresponding to operators \eqref{pamo} with $|\lambda|<1$
and $\e\ll 1$ is still almost reducible. This immediately follows from
the almost reducibility conjecture (ARC), announced by Avila in
\cite{avila}, to appear in \cite{avila1,avila2}. However, we point out
that  almost reducibility of the perturbations of the subcritical AMO
follows directly from the openness of almost reducibility and
compactness. Here we give a self-contained proof which is independent of the ARC.
\begin{Theorem}\label{sub}
Let $\alpha\in \R\backslash\Q$, $|\lambda|<1$ and $v$ be real analytic.  Then for $\e$ small enough, depending  on $\lambda,v$, the IDS of operator \eqref{pamo} is absolutely continuous.
\end{Theorem}
\begin{Remark}
In this theorem, $\alpha$ can be any irrational number, it does not need to be (strongly) Diophantine.
\end{Remark}
\subsection{Reducibility of quasiperiodic cocycles} We will only consider cocycles $(\alpha,A)$ with $\deg A=0$. A quasiperiodic $\rm C^\omega$-cocycle $(\alpha,A)$  is called {\it $\rm C^\omega$-rotations reducible} if there exists $B\in \rm C^\omega(\T,\SL(2,\R))$ and $\psi\in \rm C^\omega(\T,\R)$  such that
\begin{equation}\label{eq_reducible_def}
B(x+\alpha)^{-1}A(x)B(x)=R_{\psi(x)}.
\end{equation}
We will call a cocycle {\it almost reducible} if there exist a sequence of $B_n\in \rm C^\omega(\T,\PSL(2,\R))$ and $A_n\in SL(2,\R)$  such that
\begin{equation}
B_n(x+\alpha)^{-1}A(x)B_n(x)=A_n,
\end{equation}
moreover, $A_n\rightarrow A$ for some $A\in SL(2,\R)$.

Now we consider quasiperiodic Schr\"odinger cocycles. Define the following subset:
$$
\mathcal{R}_{\alpha,V}=\{E\in \R\colon (\alpha,S_{E}^V)\text{ is $\rm C^\omega$-rotations reducible}\}.
$$
We will sometimes drop the indices and simply use $\mathcal{R}$, if the values of the indices are clear from the context. We mention that  (almost) reducibility has been proved to be very fruitful \cite{ds,aj1,e,ak, ayz,gyzh1, ayz1,lyzz,gyz,gy,gk,houyou}.

For every $\tau>1$ and $\gamma>0$, we define
$$
\Theta^\tau_\gamma=\left\{\theta\in\T|\|2\theta+k\alpha\|_{\R/\Z}\geq \frac{\gamma}{(|k|+1)^\tau},k\in\Z\right\},
$$
$$
\Theta=\cup_{\gamma>0,\tau>1}\Theta^\tau_\gamma.
$$
Note that $\Theta$ is a subset of $\T$ of full Lebesgue measure.

\subsection{Proof of Theorem \ref{sub}} We first prove a theorem related to rotations reducibility.
\begin{Theorem}\label{rotation}
Let $\alpha\in \R\backslash\Q$, $|\lambda|<1$ and $v$ be real analytic.  Then for $\e$ small enough, depending  on $\lambda,v$, $(\alpha,S_E^{2\lambda\cos+\e v})$ is rotations reducible if $\rho(\alpha,S_E^{2\lambda\cos+\e v})\in \Theta$.
\end{Theorem}
\begin{pf}
We denote by
$$
\beta=\beta(\alpha)=\limsup\limits_{k\rightarrow \infty}-\frac{\ln\|k\alpha\|_{\R/\Z}}{|k|}.
$$
Note that if $\alpha\in\R\backslash\Q$, $(\alpha,S_E^{2\lambda\cos})$ is almost reducible for all $E\in\R$ (see Theorem 1.4 of \cite{aj1} for the case $\beta=0$ and Theorem 1.1 of \cite{avila1} for the case $\beta>0$). We then need the following two known results,
\begin{Theorem}[Corollary 1.3 of \cite{avila1}]\label{open}
Almost reducibility is stable, in the sense that it defines an open set in $\R\backslash\Q\times C^\omega(\T,SL(2,\R))$.
\end{Theorem}
Note that for $E\in \Sigma_{\lambda,\alpha}^\e$, there is $h_0>0$ such that
$$
\left|S_E^{2\lambda\cos+\e v}-S_E^{2\lambda\cos}\right|_{h_0}\leq |\e v|_{h_0}.
$$
$\Sigma_{\lambda,\alpha}^\e$ is compact and by Theorem \ref{open}, for $\e$ sufficiently small, we have that $(\alpha,S_E^{2\lambda\cos+\e v})$ is almost reducible if  $E\in \Sigma_{\lambda,\alpha}^\e$. Finally, Theorem \ref{rotation} follows from the following theorem.
\begin{Theorem}[Corollary 1.3 of \cite{avila1}]
If $(\alpha,A)$ is almost reducible and $\rho(\alpha,A)\in \Theta$, then $(\alpha,A)$ is rotations reducible.
\end{Theorem}

\end{pf}

\noindent {\bf Proof of Theorem \ref{sub}:} Let
$$
E\in \mathcal{E}=\left\{E\in\Sigma_{\lambda,\alpha}^\e:\rho(\alpha,S_E^{2\lambda\cos+\e v})\in \Theta\right\}.
$$
By Theorem \ref{rotation}, we have $(\alpha,S_E^{2\lambda\cos+\e v})$
is rotations reducible for all $E\in \mathcal{E}$. Thus by the
argument in Avila-Fayad-Krikorian \cite{afk}, $N(E)=1-2\rho(E)$ is
Lipschitz on $\mathcal{E}$. Thus the image of non-Lipschitz $E's$
under $N$ is of zero Lebesgue measure, therefore the IDS of operator \eqref{pamo} is absolutely continuous under the assumption that $|\lambda|<1$ and $\e$ sufficiently small. \footnote{We mention that this argument was first used in \cite{avila1}.}\qed

\section*{Acknowledgement}
This work was partially supported by NSF DMS-1901462, DMS-2052899, and
Simons 681675. L. Ge was also partially supported by AMS-Simons Travel Grant 2020-2022.

\end{document}